\documentclass{article}
\usepackage{amssymb}
\usepackage{graphicx}
\usepackage{graphicx}
\usepackage{float,epsfig}

\newtheorem{theorem}{Theorem}[section]
\newtheorem{lemma}[theorem]{Lemma}
\newtheorem{corollary}[theorem]{Corollary}

\newtheorem{remark}[theorem]{Remark}

\def\Capacity{\mathop{\rm Cap}\nolimits}
\def\sn{\mathop{\rm sn}\nolimits}

\long\def\symbolfootnote[#1]#2{\begingroup%
\def\thefootnote{\fnsymbol{footnote}}\footnote[#1]{#2}\endgroup}

\begin{document}
\title{\textbf{RIEMANN--SCHWARZ \\ REFLECTION PRINCIPLE AND\\  ASYMPTOTICS OF
MODULES\\ OF RECTANGULAR FRAMES}}
\author{S.R.~Nasyrov}
\date{23 of May, 2013}
\maketitle

\hskip 30cm \symbolfootnote[0]{\textsc{2010 AMS Subject
Classification.} Primary 30C20. \\
$\phantom{mm}$\textsc{Key words and phrases.} Conformal module,
doubly-connected plane domain, polygon, quadrilateral, reflection
principle, elliptic integrals, anharmonic ratio,
conver\-gen\-ce to a kernel, prime ends, asymptotics. \\
$\phantom{mm}$ The work is supported by the grants No.~12-01-97015
and No.~11-01-00762 of the Russian Fund of Basic Research.}

\begin{abstract}
We investigate asymptotical behavior of the conformal module of  a
doubly-connected domain which is the  difference of two
homo\-the\-tic re\-c\-tan\-gles under stretching it along the
abscissa axis. Thereby, we give the answer to a question put by
Prof. M.~Vuorinen.
\end{abstract}

\section{Introduction}\label{1}

In recent years an interest increases to investigation of conformal
modules of plane quadrilaterals, doubly-connected domains and
capacities of condensers with polygonal boundaries. Special
consideration is given to studying of be\-ha\-vi\-or of the modules
under various deformations of domains, their numerical
cal\-cu\-la\-tion, and asymptotics at degeneration. In this regard
we can note the review by R.~K\"uhnau~\cite{kuhnau} and the
papers~\cite{betsakos}--\cite{hakula}.\par

At first, we recall some definitions. Consider a plane
doubly-connected domain $D$ with nondegenerated boundary components.
One of its important characteristics is the conformal module $m(D)$.
There are several equivalent definitions of $m(D)$; we give some of
them.\par

If $D$ is conformally equivalent to an annulus $\{r_1<|z|<r_2\}$,
then $$
m(D):=\frac{1}{2\pi}\ln\frac{r_2}{r_1}.
$$
On the other hand,
$$ m(D):=\lambda(\Gamma),
$$
where $\lambda(\Gamma)$ is the extremal length of curve-family
$\Gamma$ consisting of all curves joining in $D$ its boundary
components. Besides,
$$
m(D):=1/\lambda(\Gamma'),
$$
$\Gamma'$ being the family of all curves in $D$ separating its
boundary components. At last,
$$
m(D):=1/\Capacity(C),
$$
where $\Capacity(C)$ is the conformal capacity of the condenser $C$
defined by~$D$.\par

Significant property of module is its invariance under conformal
mappings.  It is quasiinvariant under quasiconformal maps (see,
e.~g., \cite{ahlfors}): if $f$ is an $H$-quasicon\-for\-mal mapping
of $D$ onto $\widetilde{D}$, then
$$
\frac{1}{H}\,m(D)\le m(\widetilde{D})\le H\,m(D).
$$\par

One of the simplest $H$-quasiconformal mappings is the stretching
along the abscissa axis $f_H: x+iy\mapsto Hx+iy$. M.~Vuorinen states
the following problem: to investigate how the module $m(D)$ is
deformed under $f_H$ for sufficiently large~$H$. In particular,
which is asymptotical behavior of $m(D)$ if $D$ is the difference of
two homothetic squares?\par

The main result of the paper is\par

\begin{theorem}\label{main}
If $D_1=D_1^\sigma:=[-1,1]^2\setminus [-\sigma,\sigma]^2$, $\sigma
\in (0,1)$, $D_H=D_H^\sigma:=f_H(D_1)$, then
\begin{equation}\label{modalpha}
m(D_H^\sigma)\sim \frac{1-\sigma}{4\sigma H},\quad H\to\infty.
\end{equation}
\end{theorem}

In Section~\ref{2} we give a solution to the problem for
$\sigma=1/2$, in addition, we deduce an explicit formula for
$m(D_H)$ via elliptic integrals. It should be noted that when $H=1$
an explicit formula for $m(D_H)$ is well-known, see
Remark~\ref{bowman} below. In Section~\ref{4} the general case is
considered.  The results of Sections~\ref{2} and \ref{4} were
announced in \cite{borisova} and \cite{nasyrov}. Preliminarily, in
Section~\ref{3} we establish continuity of module of quadrilateral
under kernel convergence in the sense of Carath\'eodory.\par

To prove our results we need to recall some more definitions and
facts. A \textit{quadrilateral} is a simply-connected domain $D$,
conformally equivalent to a disk, with four marked distinct points
(prime ends) $A_k$, $1\le k\le 4$, on the boundary of $D$;
increasing of $k$ corresponds to the order in which the points occur
when we bypass $\partial D$ in the positive direction. We denote the
quadrilateral as $D(A_1,A_2,A_3,A_4)$ or simply $D$ if it is clear
which points $A_k$ are fixed. The parts of $\partial D$ lying
between $A_1$  and $A_2$, $A_3$ and $A_4$ we call \textit{horizontal
sides} of $D$, the other two parts of the boundary are
\textit{vertical sides}. Let us map conformally $D(A_1,A_2,A_3,A_4)$
onto a rectangle $[0,a]\times [0,b]$ so that the horizontal sides
are mapped onto the horizontal sides of the rectangular. The number
$$m(D):=\frac{a}{b}$$ is called the \textit{module} of $D$. It is
known (see, e.~g., \cite{ahlfors}) that $m(D)$ is equal to the
extremal length $\lambda(\Gamma)$ of the family $\Gamma$ consisting
of curves in $D$ joining its vertical sides; therefore, it is
invariant under conformal mappings and quasiinvariant under
quasiconformal ones. Besides, $m(D)=1/\lambda(\Gamma')$ where
$\Gamma'$ is the family of curves in $D$ joining its horizontal
sides.\par

Let $E:=\{|z|<1\}$, $U:=\{\Im z
>0\}$, $E^+:=E\cap U\}$, $S_{\gamma\delta}:=\{e^{i\varphi}\mid
\gamma<\varphi<\delta\}$, $0<\delta-\gamma<2\pi$. We denote by
$[a,b]$ the segment with endpoints $a$, $b\in \mathbb{C}$.\par

The elliptic integral of the first kind
$$
K(r):=\int_0^1 \frac{d\xi}{\sqrt{(1-\xi^2)(1-r^2\xi^2)}}\,.
$$
It is known (see, e.~g., \cite{avv}, \cite{akhiezer}) that
\begin{equation}\label{ell-as1}
\lim_{r\to0}\left(K(r')-\ln\frac{4}{r}\right)=0,
\end{equation}
where as usual $r'=\sqrt{1-r^2}$. From (\ref{ell-as1}) it follows
that $K(r)\sim\ln\frac{4}{r'}$ as $r\to 1$. Therefore,
\begin{equation}\label{ell-as2}
K(r)\sim\frac{1}{2}\,\ln\frac{1}{1-r}, \quad
\frac{K(r)}{K(r')}\sim\frac{1}{\pi}\,\ln\frac{1}{1-r},\quad r\to1.
\end{equation}
From (\ref{ell-as1}) we also obtain that
\begin{equation}\label{ell-as3}
\frac{K(r')}{K(r)}\sim\frac{2}{\pi}\,\ln\frac{1}{r},\quad r\to0.
\end{equation}

\begin{remark} The ring domain consisting of the unit disk minus a radial
slit from $0$ to $r$, $0<r<1$, is usually called the Gr\"otzsch ring
and its modulus is often denoted $$\mu(r) = \frac{\pi}{2}
\frac{K(r')}{K(r)}\,,$$  see~\cite{kuhnau}. The asymptotic formula
(\ref{ell-as3}) can be refined by use of the results from
\cite{avv}, Theorem~5.13.\end{remark}

\begin{remark}\label{bowman} When $H=1$ an explicit formula for $m(D^\sigma_H)$ is well-known
 (see, e.~g., \cite{bow}):
$$
m(D^\sigma_1)=
\mu\left(\left(\frac{l-l'}{l+l'}\right)^2\right),\quad
l=\mu^{-1}\left(\frac{2}{\pi}\,\frac{1-\sigma}{1+\sigma}\right),\quad
l'=\sqrt{1-l^2}.
$$
\end{remark}

\section{The case $\sigma=1/2$}\label{2} Consider the part  $Q_H$ of
$D_H$ lying in the first quarter of the plane. It is the union of
three rectangles of the same size. Let us map conformally one of the
rectangles with vertices at the points $(H+i)/2$, $H+i/2$, $H/2+i$,
and $H+i$ onto the quarter of the unit disk $U_1:=\{z\mid |z|<1, \Re
z>0, \Im z>0\}$ by the mapping $f$ so that $f((H+i)/2)=0$,
$f(H+i/2)=1$, and $f(H/2+i)=i$. Let $e^{i\kappa}=f(H+i)$.\par

By the Riemann-Schwarz reflection principle $f$ could be extended up
to the conformal mapping of the rectangle $[0,H]\times [0,1]$ onto
the unit disk $E$; we will designate the extension also through $f$.
Then $f$ maps conformally $Q_H$ onto the domain which is three
quarters of the unit disk, besides, $f(H/2)=-i$,
$f(H)=e^{-i\kappa}$, $f(i/2)=-1$, and $f(i)=-e^{-i\kappa}$.\par

\vskip 1 cm \hskip -1.4 cm \unitlength 1mm \linethickness{0.4pt}
\ifx\plotpoint\undefined\newsavebox{\plotpoint}\fi
\begin{picture}(124.25,61.306)(0,0)
\put(7.75,20.25){\framebox(52.25,40.75)[cc]{}}
\put(19.5,30.25){\framebox(28,20)[cc]{}}
\put(7.75,40.25){\line(1,0){52.25}}
\put(122.603,40.5){\line(0,1){.9106}}
\put(122.581,41.411){\line(0,1){.9087}}
\put(122.518,42.319){\line(0,1){.9047}}
\multiput(122.412,43.224)(-.0295,.179771){5}{\line(0,1){.179771}}
\multiput(122.265,44.123)(-.031516,.148505){6}{\line(0,1){.148505}}
\multiput(122.076,45.014)(-.032898,.125898){7}{\line(0,1){.125898}}
\multiput(121.845,45.895)(-.030109,.096626){9}{\line(0,1){.096626}}
\multiput(121.574,46.765)(-.031109,.085611){10}{\line(0,1){.085611}}
\multiput(121.263,47.621)(-.031865,.07643){11}{\line(0,1){.07643}}
\multiput(120.913,48.462)(-.032433,.068629){12}{\line(0,1){.068629}}
\multiput(120.524,49.285)(-.0328488,.0618902){13}{\line(0,1){.0618902}}
\multiput(120.097,50.09)(-.0331393,.0559905){14}{\line(0,1){.0559905}}
\multiput(119.633,50.874)(-.0333242,.0507645){15}{\line(0,1){.0507645}}
\multiput(119.133,51.635)(-.0334186,.046089){16}{\line(0,1){.046089}}
\multiput(118.598,52.373)(-.0334339,.0418699){17}{\line(0,1){.0418699}}
\multiput(118.03,53.084)(-.0333794,.0380343){18}{\line(0,1){.0380343}}
\multiput(117.429,53.769)(-.0332623,.0345245){19}{\line(0,1){.0345245}}
\multiput(116.797,54.425)(-.0348303,.0329421){19}{\line(-1,0){.0348303}}
\multiput(116.135,55.051)(-.0383409,.0330267){18}{\line(-1,0){.0383409}}
\multiput(115.445,55.645)(-.0421769,.0330458){17}{\line(-1,0){.0421769}}
\multiput(114.728,56.207)(-.0463957,.0329915){16}{\line(-1,0){.0463957}}
\multiput(113.986,56.735)(-.0510701,.0328539){15}{\line(-1,0){.0510701}}
\multiput(113.22,57.228)(-.0562941,.0326207){14}{\line(-1,0){.0562941}}
\multiput(112.431,57.684)(-.062191,.0322758){13}{\line(-1,0){.062191}}
\multiput(111.623,58.104)(-.068925,.031798){12}{\line(-1,0){.068925}}
\multiput(110.796,58.486)(-.076721,.031158){11}{\line(-1,0){.076721}}
\multiput(109.952,58.828)(-.095438,.033685){9}{\line(-1,0){.095438}}
\multiput(109.093,59.132)(-.109013,.032867){8}{\line(-1,0){.109013}}
\multiput(108.221,59.394)(-.126196,.031734){7}{\line(-1,0){.126196}}
\multiput(107.338,59.617)(-.14879,.030144){6}{\line(-1,0){.14879}}
\multiput(106.445,59.797)(-.180036,.027839){5}{\line(-1,0){.180036}}
\put(105.545,59.937){\line(-1,0){.9057}}
\put(104.639,60.034){\line(-1,0){.9092}}
\put(103.73,60.089){\line(-1,0){.9108}}
\put(102.819,60.102){\line(-1,0){.9104}}
\put(101.909,60.072){\line(-1,0){.908}}
\put(101.001,60){\line(-1,0){.9037}}
\multiput(100.097,59.886)(-.179491,-.03116){5}{\line(-1,0){.179491}}
\multiput(99.199,59.731)(-.148208,-.032887){6}{\line(-1,0){.148208}}
\multiput(98.31,59.533)(-.10989,-.029802){8}{\line(-1,0){.10989}}
\multiput(97.431,59.295)(-.096344,-.031){9}{\line(-1,0){.096344}}
\multiput(96.564,59.016)(-.08532,-.031898){10}{\line(-1,0){.08532}}
\multiput(95.711,58.697)(-.076133,-.03257){11}{\line(-1,0){.076133}}
\multiput(94.873,58.339)(-.068326,-.033066){12}{\line(-1,0){.068326}}
\multiput(94.053,57.942)(-.0615842,-.033419){13}{\line(-1,0){.0615842}}
\multiput(93.253,57.507)(-.055682,-.033655){14}{\line(-1,0){.055682}}
\multiput(92.473,57.036)(-.0473011,-.0316797){16}{\line(-1,0){.0473011}}
\multiput(91.716,56.529)(-.0430855,-.0318521){17}{\line(-1,0){.0430855}}
\multiput(90.984,55.988)(-.0392505,-.0319404){18}{\line(-1,0){.0392505}}
\multiput(90.277,55.413)(-.0377243,-.0337293){18}{\line(-1,0){.0377243}}
\multiput(89.598,54.806)(-.0342158,-.0335798){19}{\line(-1,0){.0342158}}
\multiput(88.948,54.168)(-.032619,-.035133){19}{\line(0,-1){.035133}}
\multiput(88.328,53.5)(-.0326712,-.0386443){18}{\line(0,-1){.0386443}}
\multiput(87.74,52.805)(-.0326548,-.0424804){17}{\line(0,-1){.0424804}}
\multiput(87.185,52.082)(-.0325616,-.0466984){16}{\line(0,-1){.0466984}}
\multiput(86.664,51.335)(-.0323808,-.0513714){15}{\line(0,-1){.0513714}}
\multiput(86.179,50.565)(-.0320994,-.056593){14}{\line(0,-1){.056593}}
\multiput(85.729,49.772)(-.0317,-.0624864){13}{\line(0,-1){.0624864}}
\multiput(85.317,48.96)(-.03116,-.069216){12}{\line(0,-1){.069216}}
\multiput(84.943,48.129)(-.033493,-.084706){10}{\line(0,-1){.084706}}
\multiput(84.608,47.282)(-.032802,-.095745){9}{\line(0,-1){.095745}}
\multiput(84.313,46.421)(-.031859,-.109312){8}{\line(0,-1){.109312}}
\multiput(84.058,45.546)(-.030567,-.126484){7}{\line(0,-1){.126484}}
\multiput(83.844,44.661)(-.028768,-.149062){6}{\line(0,-1){.149062}}
\multiput(83.672,43.766)(-.03272,-.22536){4}{\line(0,-1){.22536}}
\put(83.541,42.865){\line(0,-1){.9065}}
\put(83.452,41.958){\line(0,-1){1.8205}}
\put(83.401,40.138){\line(0,-1){.9101}}
\put(83.439,39.228){\line(0,-1){.9073}}
\multiput(83.519,38.32)(.03057,-.22566){4}{\line(0,-1){.22566}}
\multiput(83.641,37.418)(.032816,-.179195){5}{\line(0,-1){.179195}}
\multiput(83.805,36.522)(.029361,-.12677){7}{\line(0,-1){.12677}}
\multiput(84.011,35.634)(.030816,-.10961){8}{\line(0,-1){.10961}}
\multiput(84.257,34.758)(.031888,-.096054){9}{\line(0,-1){.096054}}
\multiput(84.544,33.893)(.032685,-.085022){10}{\line(0,-1){.085022}}
\multiput(84.871,33.043)(.033272,-.075829){11}{\line(0,-1){.075829}}
\multiput(85.237,32.209)(.033695,-.068018){12}{\line(0,-1){.068018}}
\multiput(85.642,31.393)(.0315588,-.0568963){14}{\line(0,-1){.0568963}}
\multiput(86.083,30.596)(.03189,-.0516775){15}{\line(0,-1){.0516775}}
\multiput(86.562,29.821)(.0321152,-.0470065){16}{\line(0,-1){.0470065}}
\multiput(87.076,29.069)(.0322487,-.0427895){17}{\line(0,-1){.0427895}}
\multiput(87.624,28.341)(.0323016,-.0389538){18}{\line(0,-1){.0389538}}
\multiput(88.205,27.64)(.0322828,-.0354422){19}{\line(0,-1){.0354422}}
\multiput(88.819,26.967)(.0321997,-.032209){20}{\line(0,-1){.032209}}
\multiput(89.463,26.323)(.0354328,-.0322931){19}{\line(1,0){.0354328}}
\multiput(90.136,25.709)(.0389444,-.0323129){18}{\line(1,0){.0389444}}
\multiput(90.837,25.127)(.0427802,-.0322611){17}{\line(1,0){.0427802}}
\multiput(91.564,24.579)(.0469972,-.0321289){16}{\line(1,0){.0469972}}
\multiput(92.316,24.065)(.0516683,-.031905){15}{\line(1,0){.0516683}}
\multiput(93.091,23.586)(.0568871,-.0315753){14}{\line(1,0){.0568871}}
\multiput(93.888,23.144)(.068008,-.033715){12}{\line(1,0){.068008}}
\multiput(94.704,22.74)(.075819,-.033294){11}{\line(1,0){.075819}}
\multiput(95.538,22.373)(.085012,-.032709){10}{\line(1,0){.085012}}
\multiput(96.388,22.046)(.096044,-.031916){9}{\line(1,0){.096044}}
\multiput(97.252,21.759)(.109601,-.030848){8}{\line(1,0){.109601}}
\multiput(98.129,21.512)(.126761,-.029398){7}{\line(1,0){.126761}}
\multiput(99.016,21.307)(.179186,-.032868){5}{\line(1,0){.179186}}
\multiput(99.912,21.142)(.22565,-.03064){4}{\line(1,0){.22565}}
\put(100.815,21.02){\line(1,0){.9073}}
\put(101.722,20.939){\line(1,0){.9101}}
\put(102.632,20.901){\line(1,0){.9109}}
\put(103.543,20.905){\line(1,0){.9097}}
\put(104.453,20.951){\line(1,0){.9066}}
\multiput(105.359,21.04)(.22537,.03265){4}{\line(1,0){.22537}}
\multiput(106.261,21.171)(.14907,.028725){6}{\line(1,0){.14907}}
\multiput(107.155,21.343)(.126493,.03053){7}{\line(1,0){.126493}}
\multiput(108.041,21.557)(.109321,.031827){8}{\line(1,0){.109321}}
\multiput(108.915,21.811)(.095755,.032774){9}{\line(1,0){.095755}}
\multiput(109.777,22.106)(.084716,.033468){10}{\line(1,0){.084716}}
\multiput(110.624,22.441)(.069225,.03114){12}{\line(1,0){.069225}}
\multiput(111.455,22.815)(.0624956,.0316819){13}{\line(1,0){.0624956}}
\multiput(112.267,23.226)(.0566024,.032083){14}{\line(1,0){.0566024}}
\multiput(113.06,23.676)(.0513808,.0323659){15}{\line(1,0){.0513808}}
\multiput(113.83,24.161)(.0467079,.032548){16}{\line(1,0){.0467079}}
\multiput(114.578,24.682)(.0424898,.0326425){17}{\line(1,0){.0424898}}
\multiput(115.3,25.237)(.0386538,.03266){18}{\line(1,0){.0386538}}
\multiput(115.996,25.825)(.0351425,.0326088){19}{\line(1,0){.0351425}}
\multiput(116.664,26.444)(.0335897,.0342061){19}{\line(0,1){.0342061}}
\multiput(117.302,27.094)(.0319644,.0357296){19}{\line(0,1){.0357296}}
\multiput(117.909,27.773)(.0319518,.0392412){18}{\line(0,1){.0392412}}
\multiput(118.484,28.479)(.0318646,.0430763){17}{\line(0,1){.0430763}}
\multiput(119.026,29.212)(.0316934,.0472919){16}{\line(0,1){.0472919}}
\multiput(119.533,29.968)(.0336712,.0556722){14}{\line(0,1){.0556722}}
\multiput(120.004,30.748)(.0334369,.0615745){13}{\line(0,1){.0615745}}
\multiput(120.439,31.548)(.033085,.068316){12}{\line(0,1){.068316}}
\multiput(120.836,32.368)(.032592,.076123){11}{\line(0,1){.076123}}
\multiput(121.195,33.205)(.031923,.085311){10}{\line(0,1){.085311}}
\multiput(121.514,34.058)(.031028,.096335){9}{\line(0,1){.096335}}
\multiput(121.793,34.926)(.029834,.109881){8}{\line(0,1){.109881}}
\multiput(122.032,35.805)(.03293,.148198){6}{\line(0,1){.148198}}
\multiput(122.229,36.694)(.031212,.179482){5}{\line(0,1){.179482}}
\put(122.385,37.591){\line(0,1){.9037}}
\put(122.5,38.495){\line(0,1){.908}}
\put(122.572,39.403){\line(0,1){1.0971}}
\put(102.68,59.93){\line(0,-1){.95}}
\put(102.68,58.03){\line(0,-1){.95}}
\put(102.68,56.13){\line(0,-1){.95}}
\put(102.68,54.23){\line(0,-1){.95}}
\put(102.68,52.33){\line(0,-1){.95}}
\put(102.68,50.43){\line(0,-1){.95}}
\put(102.68,48.53){\line(0,-1){.95}}
\put(102.68,46.63){\line(0,-1){.95}}
\put(102.68,44.73){\line(0,-1){.95}}
\put(102.68,42.83){\line(0,-1){.95}}
\put(102.68,40.93){\line(1,0){.9875}}
\put(104.655,40.93){\line(1,0){.9875}}
\put(106.63,40.93){\line(1,0){.9875}}
\put(108.605,40.93){\line(1,0){.9875}}
\put(110.58,40.93){\line(1,0){.9875}}
\put(112.555,40.93){\line(1,0){.9875}}
\put(114.53,40.93){\line(1,0){.9875}}
\put(116.505,40.93){\line(1,0){.9875}}
\put(118.48,40.93){\line(1,0){.9875}}
\put(120.455,40.93){\line(1,0){.9875}}
\put(103,40.75){\line(0,-1){19.25}} \put(103,41){\line(-1,0){20}}
\put(57.25,3){\line(0,1){.25}}
\put(63.75,9){\makebox(0,0)[cc]{Fig.~1}}
\put(120.25,50){\circle*{1}}
\put(124.25,51.5){\makebox(0,0)[cc]{$e^{i\kappa}$}}
\put(120,30.5){\circle*{1}}
\put(124.25,29.75){\makebox(0,0)[cc]{$e^{-i\kappa}$}}
\put(86,50){\circle*{1}}
\put(81.25,51){\makebox(0,0)[cc]{$-e^{-i\kappa}$}}
\put(33.25,61){\line(0,-1){40.75}} \put(103,41){\line(0,-1){20}}
\put(47.335,50.278){\line(0,1){.9651}}
\put(47.335,52.209){\line(0,1){.9651}}
\put(47.335,54.139){\line(0,1){.9651}}
\put(47.335,56.069){\line(0,1){.9651}}
\put(47.335,57.999){\line(0,1){.9651}}
\put(47.335,59.93){\line(0,1){.9651}}
\put(47.335,50.068){\line(1,0){.9622}}
\put(49.26,50.068){\line(1,0){.9622}}
\put(51.184,50.068){\line(1,0){.9622}}
\put(53.108,50.068){\line(1,0){.9622}}
\put(55.033,50.068){\line(1,0){.9622}}
\put(56.957,50.068){\line(1,0){.9622}}
\put(58.881,50.068){\line(1,0){.9622}}
\put(47.511,50.138){\circle*{1.}} \put(47.511,40.258){\circle*{1.}}
\put(59.914,40.153){\circle*{1.}} \put(60.019,60.86){\circle*{1.}}
\put(33.215,50.244){\circle*{1.}} \put(33.321,60.965){\circle*{1.}}
\put(103.01,41.099){\circle*{1.}} \put(103.01,21.022){\circle*{1.}}
\put(83.354,40.994){\circle*{1.}}
\put(57.917,43.){\makebox(0,0)[cc]{$H$}}
\put(43.503,43.){\makebox(0,0)[cc]{$H/2$}}
\put(32.093,43.){\makebox(0,0)[cc]{$0$}}
\put(101.093,43.){\makebox(0,0)[cc]{$0$}}
\put(101.093,18.){\makebox(0,0)[cc]{$-i$}}
\put(103.093,63.){\makebox(0,0)[cc]{$i$}}
\put(80.093,42.){\makebox(0,0)[cc]{$-1$}}
\put(125.093,42.){\makebox(0,0)[cc]{$1$}}
\put(55.00,58.309){\makebox(0,0)[cc]{$H+i$}}
\put(35.761,58.704){\makebox(0,0)[cc]{$i$}}
\put(36.267,52.61){\makebox(0,0)[cc]{$i/2$}}
\end{picture}
\vskip -0.3 cm

The function $g(\zeta):=(if(\zeta))^{2/3}$ maps conformally $Q_H$
onto $E^+$ at an ap\-pro\-pri\-ate choice of a branch of the power
function. We have $g(H/2)=1$, $g(H)=e^{i\alpha}$, $f(i/2)=-1$,
$f(i)=-e^{-i\beta}$ where
\begin{equation}\label{angles}
\alpha=(\pi-2\kappa)/3,\quad \beta=\pi-2\kappa/3.
\end{equation}\par

The module of $D_H$, by the symmetry principle for quasiconformal
mappings (see, e.~g., \cite{ahlfors}) is equal
\begin{equation}\label{quater}
m(D_H)={1}/{(4\lambda(\Gamma))}
\end{equation}
where $\Gamma$ is the family of all curves in $Q_H$ which join
$[H/2,H]$ and $[i/2,i]$. Because of conformal invariance of the
module we obtain
\begin{equation}\label{prime}
\lambda(\Gamma)=\lambda(\Gamma'),\end{equation} where $\Gamma'$ is a
family of all curves in $E^+$ connecting $S_{0\alpha}$ to
$S_{\beta\pi}$. By  the symmetry principle,
\begin{equation}\label{sym}
\lambda(\Gamma')=2\lambda(\widetilde{\Gamma}) \end{equation} where
$\widetilde{\Gamma}$ is the family of all curves which join
$S_{-\alpha,\alpha}$ and $S_{\beta,2\pi-\beta}$ in $E^+$.\par

Let us map conformally the unit disk $E$ onto the upper half-plane
$U$ so that the points $e^{i\beta}$, $e^{-i\beta}$, $e^{-i\alpha}$,
and $e^{i\alpha}$ are mapped on $-1/l$, $-1$, $1$, and $1/l$, $l>1$.
Now we express $l$ through $\alpha$ and $\beta$. Equating the
anharmonic ratios
$$ \frac{-1/l+1}{-1/l-1}\cdot
\frac{1/l-1}{1/l+1}=\frac{e^{i\beta}-e^{-i\beta}}
{e^{i\beta}-e^{-i\alpha}}\cdot\frac{e^{i\alpha}-e^{-i\alpha}}
{e^{i\alpha}-e^{-i\beta}}
$$
we obtain
\begin{equation}\label{el}
l=\frac{\sqrt{1-\cos(\alpha+\beta)}-\sqrt{2\sin\beta\sin
\alpha}}{\sqrt{1-\cos(\alpha+\beta)}+\sqrt{2\sin\beta\sin
\alpha}}\,.
\end{equation}\par

Besides,
\begin{equation}\label{length}
\lambda(\widetilde{\Gamma})=\frac{2K(l)}{K(l')}\,.\end{equation}\par

From (\ref{quater}), (\ref{prime}), (\ref{sym}), and (\ref{length})
we have
\begin{equation}\label{modd}
m(D_H)=\frac{K(l')}{16K(l)}\,.
\end{equation}\par

Now we find the relation between $\kappa$ and $H$. For this purpose
we map conformally $E$ onto the upper half-plane $U$ by a function
$\varphi$ so that $-e^{-i\kappa}$, $-e^{i\kappa}$,
$e^{-i\kappa}$, and $e^{i\kappa}$ are mapped on $-1/k$, $-1$, $1$,
and $1/k$. We note that $k$ satisfies the condition
\begin{equation}\label{k}
\frac{2K(k)}{K(k')}=H
\end{equation}
because the quadrilateral, which is the upper half-plane $U$ with
fixed points $-1/k$, $-1$, $1$, and $1/k$, is conformally equivalent
to the rectangle of length $H$ and height $1$ under the mapping
$\varphi\circ f$. From the equality of cross-ratios
$$
\frac{-1/k+1}{-1/k-1}\cdot
\frac{1/k-1}{1/k+1}=\frac{-e^{-i\kappa}+e^{i\kappa}}{-e^{-i\kappa}-e^{i\kappa}}\cdot
\frac{e^{i\kappa}-e^{-i\kappa}}{e^{i\kappa}-e^{-i\kappa}}
$$
we have
\begin{equation}\label{kappa}
\kappa=\arcsin \frac{1-k}{1+k}\,.
\end{equation}\par

Therefore, we prove
\begin{theorem}\label{1/2} For $\sigma=1/2$ the module of $D_H=D_H^\sigma$ is defined by
(\ref{modd}) where $l$ is found from (\ref{el}) taking into account
(\ref{angles}), (\ref{k}), and (\ref{kappa}).
\end{theorem}\par

\begin{corollary} We have
$$m(D_H)\sim \frac{1}{4H},\quad H\to\infty.$$
\end{corollary} Actually, from (\ref{kappa}) it follows that $\kappa\sim (1-k)/2$ as $H\to\infty$.
Now, taking into account (\ref{el}) and (\ref{angles}), we obtain
$$
1-l\simeq\sqrt{\sin\beta}=\sqrt{\sin
\frac{2\kappa}{3}}\simeq\sqrt{1-k}. $$ Thus, using (\ref{modd}),
(\ref{k}), and (\ref{ell-as2}), we have
$$
m(D_H)= \frac{K(l')}{16K(l)}\sim \frac{\pi}{16\ln[(1-l)^{-1}]}\sim
\frac{\pi}{8\ln[(1-k)^{-1}]}\sim
 \frac{1}{4H}, \quad H\to\infty.$$

\section{Convergence of domains and their modules}\label{3}

At first, we recall some results of the theory of prime ends of a
sequence of domains converging to a kernel~\cite{suvorov}.\par

Consider a sequence of simply-connected domains $G_n$ on the Riemann
sphere $\overline{\mathbb{C}}$ converging to a kernel $G$ with
respect to a fixed point $S_0\in \overline{\mathbb{C}}$. We assume
that the boundaries of $G_n$ and $G$ are nondegenerate, i.~e., each
of them contains more than one point. Further we provide $G_n$ and
$G$ with metrics induced from the sphere; in the case when all $G_n$
are contained in a fixed Euclidean disk it is possible to change the
spherical metrics by the Euclidean one.\par

Consider a section $\gamma$ of $G$, i.~e., a Jordan arc in ${G}$
with endpoints on $\partial G$. Let $\gamma_n$ be a section of
$G_n$. We say that the sequence $(\gamma_n)$ is \textit{a section of
the sequence $\widetilde{G}:=(G_n)$ lying over} $\gamma$ if the
following conditions are fulfilled:
\begin{enumerate}
\item For any neighborhood $U$ of $\gamma$ there exists $n_0$ such
that $\gamma_n$ lies in $U$ for any $n\ge n_0$;
\item if points $p_1$, $p_2\in G$ are separated by $\gamma$ in  $D$,
then there exists $n_1$ such that  $p_1$ and $p_2$ are separated  by
$\gamma_n$ in $G_n$ for  any $n\ge n_1$.
\end{enumerate}\par

In \cite{suvorov} it is proved\par

\begin{lemma}\label{section}
For any section $\gamma$ of $G$ there exists a section
$\widetilde{\gamma}:=(\gamma_n)$ of $\widetilde{G}$ lying
over~$\gamma$.
\end{lemma}\par

Let $(\gamma_m)$ be a chain of sections of $G$ which defines a prime
end $P$ of $G$, and $\widetilde{\gamma}_m:=(\gamma_{mn})$ is a
section of $\widetilde{G}$ lying over~$\gamma_m$. It is possible to
introduce a natural relation of equivalence on the set of all such
sequences $(\widetilde{\gamma}_m)$ (in more detail see
\cite{suvorov}); the classes of equivalence $\widetilde{P}$ of such
sequences are called \textit{prime ends} of~$\widetilde{G}$.\par

If a sequence of sections $(\gamma_m)$ defines a prime end $P$ of
$G$ and prime end $\widetilde{P}$ of $\widetilde{G}$ contains
$(\widetilde{\gamma}_m)$ where $\widetilde{\gamma}_m$ is a section
of $\widetilde{G}$ lying over $\gamma_m$, then we will say that
$\widetilde{P}$ is \textit{the prime end of $\widetilde{G}$
corresponding to prime end} $P$ of $G$. The described correspondence
$\Phi:P\mapsto \widetilde{P}$ is a bijection between the set of
prime ends of $G$ and  the set of prime ends of $\widetilde{G}$.\par

Now consider a sequence $(P_n)$ where $P_n$ is a point or a boundary
prime end of $G_n$. Let $\widetilde{P}=\Phi(P)$ be the prime end of
$\widetilde{G}$ corresponding to a prime end  $P$ of $G$ and $P$ is
defined by a sequence $(\gamma_m)$. Let $\widetilde{P}$ is defined
by $(\widetilde{\gamma}_m)$, where $\widetilde{\gamma}_m$ lies over
$\gamma_m$ for any $m$, and
$\widetilde{\gamma}_m:=(\gamma_{mn})$.\par

We will say that \textit{the sequence $(P_n)$ converges to}
$\widetilde{P}$ if for any $m$ there exists $n_0$ such that
$\gamma_{mn}$ separates $P_n$ from $S_0$ in $G_m$ for any $n\ge
n_0$.\par

Let $f_n$ and $f$ be conformal mappings of $E$ onto $G_n$ and $G$,
extended to $\overline{E}$ up to homeomorphisms of prime ends. Let
$P_n$ be a sequence consisting of points or boundary prime ends of
$G_n$, and let $P$ be a boundary prime end of $G$. Denote
$\zeta_n=f_n^{-1}(P_n)$, $\zeta_0=f^{-1}(P)$.\par

Theorem~6, p.~75, from \cite{suvorov} implies\par

\begin{theorem}\label{regseq} A sequence $(P_n)$ converges to the prime
end  $\widetilde{P}$  of \ $\widetilde{G}$ cor\-res\-pon\-d\-ing to
$P$ if and only if $\zeta_n\to \zeta_0$. \end{theorem}\par

From Theorem~\ref{regseq} we deduce the following result on
continuity of the module of quadrilateral under kernel convergence
of domains.\par

\begin{theorem}\label{quadr} Let $A_n$, $B_n$, $C_n$, and $D_n$ be distinct
boundary prime ends of~$G_n$. Let the sequences $(A_n)$, $(B_n)$,
$(C_n)$, and $(D_n)$ converge to the prime ends $\widetilde{A}$,
$\widetilde{B}$, $\widetilde{C}$ and $\widetilde{D}$ of
$\widetilde{G}$ corresponding to distinct prime ends $A$, $B$, $C$,
and $D$ of~$G$. Then $m(G_n(A_n,B_n,C_n,D_n))\to
m(G(A,B,C,D))$.\end{theorem}\par

Actually, taking into account conformal invariance of the module of
quad\-ri\-la\-te\-ral, by Theorem~\ref{regseq} we have
$$m(G_n(A_n,B_n,C_n,D_n))=m(E(a_n,b_n,c_n,d_n))\to$$$$\to
E(a,b,c,d)=m(G(A,B,C,D))$$ where $a_n$, $b_n$, $c_n$, and $d_n$ are
preimages of $A_n$, $B_n$, $C_n$, and $D_n$ under the map $f_n$, and
$a$, $b$, $c$, and $d$ are preimages of $A$, $B$, $C$, and $D$ under
the map~$f$.\par

\begin{remark}\label{nonuniv} The statement of Theorem~\ref{quadr}
is valid not only for univalent domains, it
is true for $p$-valent domains (Riemann surfaces) as well (see,
e.~g.,  \cite{nasyrovm}). \end{remark}

\section{General case}\label{4} As in the case $\sigma=1/2$, consider the part $Q_H$
of $D_H$ lying in the first quarter of the plane. Let us shift $Q_H$
to the left on the value $\sigma H$; as a result we receive the
domain $\widetilde{Q}_H$. By the symmetry principle, the module of
$D_H$ is equal to
\begin{equation}\label{mod}
m(D_H)=\frac{1}{4\lambda(\Gamma)}
\end{equation}
where $\Gamma$ is the family of curves joining $[0,(1-\sigma)H]$ and
$[-\sigma H+i\sigma,-\sigma H+i]$ in  the quadrilateral
$\widetilde{Q}_H$.\par

Now consider the domain
$$\widetilde{Q}:=\cup_{H>0}\widetilde{Q}_H=(\mathbb{R}\times (0,1))\setminus
((-\infty,0]\times (0,\sigma]).$$ Let us map conformally
$\widetilde{Q}$ onto the horizontal strip $\mathbb{R}\times (0,1)$
with keeping the infinite prime ends and the origin. For this
purpose we map conformally the upper half-plane $U$ onto
$\widetilde{Q}$ and $G=\{0<\Im \omega<1\}$ by the functions
$$z=C \int_1^\zeta\sqrt{\frac{\zeta-s}{\zeta-1}}\,
\,\frac{d\zeta}{\zeta}, \quad \omega=\frac{\ln\zeta}{\pi},
$$
where $C>0$, $0<s<1$.\par

Near $\zeta=0$ we have
$$z=C\sqrt{s}\ln \zeta+\sum_{k=0}^\infty \sigma_k \zeta^k.$$
Since the intersection of the domain $\widetilde{Q}$ and the left
half-plane is a half-strip of width $(1-\sigma)/\pi$, we have
$C\sqrt{s}=(1-\sigma)/\pi$. Near $\zeta=\infty$
$$
z=C\ln \zeta+\sum_{k=0}^\infty \frac{\beta_k} {\zeta^k},
$$
therefore,  similarly we obtain $C=1/\pi$. Then
$s=(1-\sigma)^2$ and
$$
z=\frac{1}{\pi}\,
\int_1^\zeta\sqrt{\frac{\zeta-(1-\sigma)^2}{\zeta-1}}.
$$\par

Considering it we conclude that for sufficiently large $M>0$ we have
\begin{equation}\label{ser1}
z=(1-\sigma)\omega+\sum_{k=0}^\infty \sigma_k e^{k\pi\omega}.
\end{equation}
in the half-plane $\Pi_\mu^-:=\{\Re\omega<-M, \,0<\Im\omega<1\}$,
and
\begin{equation}\label{ser2} z=\omega+\sum_{k=0}^\infty
\beta_k e^{-k\pi\omega}.
\end{equation}
in the half-plane $\Pi_ \mu^+:=\{\Re\omega>M,
\,0<\Im\omega<1\}$.\par

From rectilinearity of the boundary arcs of the domains and the
Riemann-Schwarz reflection principle it follows that convergence of
the series (\ref{ser1}) and (\ref{ser2}) is uniform in the closed
half-planes $\overline{\Pi_ \mu^-}$ and $\overline{\Pi_ \mu^+}$.\par

From (\ref{ser1}) and (\ref{ser2}) we deduce that on the vertical
segments in $\widetilde{Q}$, lying on the lines $\Re
\omega=-\widetilde{\sigma}\,H$, where
$$\widetilde{\sigma}=\frac{\sigma}{1-\sigma},$$
we have
$$
\Re z(\omega)=-\widetilde{\sigma}\,H+O(1), \quad H\to\infty.
$$
In the same way, on the segments,  lying on the lines $\Re
\omega=(1-\sigma)H$,
$$
\Re z=(1-\sigma) H+O(1), \quad H\to\infty,
$$
Therefore,
\begin{equation}\label{equiv}
\lambda(\Gamma)\sim m(\widetilde{P})
\end{equation}
where $\widetilde{P}$ is the quadrilateral which is the rectangle
$$
P:=\left[-\widetilde{\sigma}\,H,(1-\sigma) H\right]\times
[0,1]
$$ with the segments
$[-\widetilde{\sigma}\,H,-\widetilde{\sigma}\,H+i]$ and
$[0,(1-\sigma) H]$ as vertical sides.\par

Let us map $U$ onto $P$ by the function $$ z=C\int_0^\zeta
\frac{d\xi}{\sqrt{(1-\xi^2)(1-k^2\xi^2)}}+C_1
$$
where
$$
C_1=\frac{(1-\sigma)^2-\sigma}{1-\sigma}\,H,\quad
C=\frac{(1-\sigma)^2+\sigma}{2 K(k)(1-\sigma)},
$$
and
$k\in (0,1)$ is defined from the relation
\begin{equation}\label{k2}
\frac{2 K(k)}{K(k')}=\frac{(1-\sigma)^2+\sigma}{1-\sigma}\, H.
\end{equation}
The mapping  takes the points $-\widetilde{\sigma}\,H+i$,
$\widetilde{\sigma}\,H$, $0$, and  $(1-\sigma)H$, i.~e., the
vertices of quadrilateral $P$, into $(-1/k)$, $(-1)$, $a$, and $1$,
where
\begin{equation}\label{sn}
a=\sn\left[\frac{\sigma-(1-\sigma)^2}{\sigma+(1-\sigma)^2}\,K(k),k\right].
\end{equation}
Here $\sn[\,\cdot\,,k]$ is the Jacobi elliptic sine corresponding to
the parameter $k$ (see, e.~g.,~\cite{akhiezer}).\par

Now we map the upper half-plane $U$ onto itself conformally so that
the points $-1/k$, $-1$, $a$, and $1$ move into $-1/\nu$, $-1$, $1$,
and $1/\nu$ ($0<\nu<1)$. Then the module of the quadrilateral
$\widetilde{P}$ is equal to
\begin{equation}\label{pmod}
m(\widetilde{P})=\frac{2 K(\nu)}{K(\nu{\,}')}
\end{equation}
where  $\nu$ is defined from the equality of the cross-ratios:
$$
\frac{1/\nu-1}{1/\nu+1}\cdot\frac{-1/\nu+1}{-1/\nu-1}=\frac{1-a}{1+1}\cdot
\frac{-1/k+1}{-1/k-a} $$ or
\begin{equation}\label{el2}
\frac{1-\nu}{1+\nu}=\sqrt{\frac{1-a}{1+ka}}\cdot\sqrt{\frac{1-k}{2}}.
\end{equation}
Therefore, for finding the asymptotics of $m(P)$  we need to know
the asymptotic behavior of $a$ as $H\to\infty$. It is possible to do
using (\ref{sn}), but we prefer to apply geometric considerations
which are based on rectilinearity of the boundary arcs and the
reflection principle. Let us prove the following auxiliary
statement.\par

\begin{lemma}\label{slit} Let $Q$ be a quadrilateral which is the square $[0,1]^2$
with vertices at the points $c$, $1$, $1+i$, and $i$, where $c\in
(0,1)$. Denote $Q_H=f_H(Q)$. Then
$$m(Q_H)\sim (1-c)H,\quad H\to\infty.$$
\end{lemma}\par

Proof. Let $\widetilde{Q}_H=(1/H)Q_H$. Since $m(\widetilde{Q}_H)$ is
a monotonic function of $H$, it is sufficient to consider the
sequence $H_n=2n$ and to prove that $$m(Q_{H_n})\sim (1-c)H_n,\quad
n\to\infty.$$ For short we will write $Q_n$ instead
of~$Q_{H_n}$.\par\vskip 0.8 cm

\hskip 0.4 cm
\unitlength 1mm
\linethickness{0.4pt}
\ifx\plotpoint\undefined\newsavebox{\plotpoint}\fi
\begin{picture}(83.262,61.96)(0,0)
\put(31.75,12.5){\framebox(42.75,39)[cc]{}}
\put(83.262,12.64){\vector(1,0){.07}}\put(24.042,12.551){\line(1,0){59.22}}
\put(55.066,12.351){\circle*{.729}}
\put(31.643,61.96){\vector(0,1){.07}}\put(31.731,6.894){\line(0,1){55.066}}
\put(74.623,12.351){\circle*{.637}}
\put(54.889,51.642){\circle*{.559}}
\put(31.731,12.463){\circle*{.53}}
\put(31.731,51.53){\circle*{.729}}
\put(55.066,13.972){\makebox(0,0)[cc]{$c$}}
\put(55.154,54.475){\makebox(0,0)[cc]{$c+i$}}
\put(74.154,54.475){\makebox(0,0)[cc]{$1+i$}}
\put(74.60,51.60){\circle*{.559}}
\put(28.726,19.18){\makebox(0,0)[cc]{$\frac{i}{n}$}}
\put(31.661,15.663){\line(1,0){1.9526}}
\put(35.566,15.663){\line(1,0){1.9526}}
\put(39.471,15.663){\line(1,0){1.9526}}
\put(43.376,15.663){\line(1,0){1.9526}}
\put(47.281,15.663){\line(1,0){1.9526}}
\put(51.187,15.663){\line(1,0){1.9526}}
\put(55.092,15.663){\line(1,0){1.9526}}
\put(58.997,15.663){\line(1,0){1.9526}}
\put(62.902,15.663){\line(1,0){1.9526}}
\put(66.807,15.663){\line(1,0){1.9526}}
\put(70.713,15.663){\line(1,0){1.9526}}
\put(29.08,51.354){\makebox(0,0)[cc]{$i$}}
\put(29.964,61.165){\makebox(0,0)[cc]{$y$}}
\put(82.555,10.43){\makebox(0,0)[cc]{$x$}}
\put(55,1){\makebox(0,0)[cc]{Fig.~2}}
\put(65,35){\makebox(0,0)[cc]{$G_n$}}
\put(54.801,44.813){\line(-1,0){23.1578}}
\put(54.978,38.184){\line(-1,0){23.1578}}
\put(54.801,32.085){\line(-1,0){23.1578}}
\put(54.801,25.191){\line(-1,0){23.1578}}
\put(54.889,19.18){\line(-1,0){23.1578}}
\put(54.801,44.283){\line(-1,0){23.069}}
\put(54.978,37.653){\line(-1,0){23.069}}
\put(54.801,31.555){\line(-1,0){23.069}}
\put(54.801,24.66){\line(-1,0){23.069}}
\put(54.889,18.65){\line(-1,0){23.069}}
\put(54.869,25.122){\line(0,-1){.473}}
\put(54.869,19.078){\line(0,-1){.473}}
\put(54.869,32.059){\line(0,-1){.473}}
\put(54.921,38.208){\line(0,-1){.473}}
\put(54.869,44.778){\line(0,-1){.473}}
\put(31.744,18.92){\circle*{.76}}
\end{picture}
\vskip 0.6 cm\par

Consider the domains $Q_n^1$, $Q_n^2,\ldots, Q_n^{2n}$, where
$$Q_n^k=[0,1]\times \left[\frac{k-1}{2n},\frac{k}{2n}\right].
$$ Let us glue
$Q_n^k$ and $Q_n^{k+1}$ along $\{(x,y) \mid 0\le x\le 1,\,
y=j/(2n)\}$ for odd $k$, and along $\{(x,y)\mid c\le x\le 1,\,
y=j/(2n)\}$ for even~$k$. As a result, we obtain the domain $G_n$
which is the unit square with $(n-1)$ horizontal slits (Fig.~2).\par

We will consider $G_n$ as a quadrilateral with vertices $c$, $1$,
$c+i$, and $i$. By the symmetry principle, $m(G_n)= m(Q_n)/(2n)$.
The domains $G_n$ converge to the rectangule $G:=[c,1]\times [0,1]$
as $n\to \infty$, and the sequences of their vertices converge to
four distinct prime ends of $\widetilde{G}=(G_n)$ corresponding to
the vertices of $G$. To show that we take  $(1/4)$ of concentric
circles with radius $r_m\to 0$ as sections $\gamma_m$ which define
prime end being a vertex of $G$. Let $\gamma_{mn}$ be a section of
$G_n$ which is the union of $\gamma_m$ and, if it is necessary, a
segment connecting one of its endpoint to one of the nearest points
of $\partial G_n$. By Theorem~\ref{quadr}, $m(G_n)\to m(G)=1-c$, and
Lemma~\ref{slit} is proved.\par

\begin{remark}\label{long} The quadrilateral $Q_H$ could be considered
as a generalized long quadri\-la\-te\-ral. Asymptotics of the
modules of long quadrilaterals was investigated in
\cite{gaier_hayman1}, \cite{gaier_hayman2}, \cite{kuhnau2},
\cite{laugesen}, \cite{falcao}, and other papers where various
methods for com\-pu\-t\-ing the modules were suggested.
\end{remark}\par

Consider the quadrilateral  $P^\ast$ which is the rectangle $P$ with
the segments $[0, (1-\sigma)H]$ and $[-\widetilde{\sigma}\,H+i,
(1-\sigma)H+i]$ as horizontal sides. Taking into account conformal
invariance of the module, by Lemma ~\ref{slit} we have
\begin{equation}\label{ast1}
m(P^\ast)\sim (1-\sigma)H, \quad H\to\infty.
\end{equation}\par

Now we can describe the behavior of $a$ as $H\to\infty$. Let us map
$P^\ast$  conformally onto $U$ such that the points $(1-\sigma)H+i$,
$-\widetilde{\sigma}H+i$, $0$, and $(1-\sigma)H$ are mapped into
$-1/\mu$, $-1$, $1$, and $1/\mu$, $0<\mu<1$. Then
\begin{equation}\label{ast2}
m(P^\ast)=\frac{K(\mu')}{2 K(\mu)},
\end{equation}
We should note that because of $m(P^\ast)\to\infty$ as $H\to\infty$,
by (\ref{ast2}) and (\ref{ell-as3}), we have $\mu\to0$,
$H\to\infty$. Comparing the anharmonic ratios of the points  in
$\partial U$, corresponding to each other under conformal
automorphism, we obtain $$
\frac{1-1/\mu}{1+1}\cdot\frac{-1/\mu+1}{-1/\mu-1/\mu}=\frac{a-1}{a+1/k}\cdot
\frac{1/k+1/k}{1/k-1}
$$ or
$$ \frac{1-a}{1+ak}\cdot\frac{2k}{1-k}=\frac{(\mu-1)^2}{4\mu}.
$$
Therefore, taking into account that $a$, $k\to1$ as $H\to\infty$ we
have
\begin{equation}\label{a}
1-a\sim \frac{1-k}{\mu}.
\end{equation}
By (\ref{ast1}), (\ref{ast2}), and (\ref{ell-as2}),
\begin{equation}\label{mu}
\frac{1}{(1-\sigma)H}\sim\frac{2K(\mu)}{K(\mu')}\sim\frac{4}{\pi}\ln
\frac{1}{\mu}\,.
\end{equation}\par

Now we use the asymptotic behavior ~(\ref{ell-as2}) of the elliptic
integrals. With use of that and by (\ref{k2}) $$ \ln
\frac{1}{1-k}\sim \pi\,\frac{K(k)}{K(k')} \sim\frac{\pi}{2}\cdot
\frac{(1-\sigma)^2+\sigma}{1-\sigma}\, H.
$$
Now from
(\ref{ell-as2}), (\ref{pmod}), (\ref{el2}), (\ref{a}), and
(\ref{mu}) we obtain
$$
m(\widetilde{P})=\frac{2K(\nu)}{K(\nu')}\sim \frac{2}{\pi}\ln
\frac{1}{1-\nu}\sim \frac{2}{\pi}\ln
\frac{1+ak}{1-a}+\frac{1}{\pi}\ln \frac{1}{1-k}$$
$$\sim
\frac{2}{\pi}\ln \frac{1}{1-k}-\frac{1}{\pi}\ln \frac{1}{\mu} \sim
(1-\sigma+\widetilde{\sigma})H- (1-\sigma)H=\widetilde{\sigma}H. $$
Because of (\ref{mod}) and (\ref{equiv}) it completes the proof
of~(\ref{modalpha}).\vskip 0.3 cm\par

\textit{Acknowledgement.} The author expresses his gratitude to
Prof. M.~Vuorinen for attraction attention to the problem and useful
comments, and to E.~V.~Bo\-ri\-so\-va for help in checking the
argument in the proof of Theorem~\ref{1/2}.

\vskip 0.3 cm \textsc{Institute of mathematics and mechanics, Kazan
(Volga Region) Federal University, e-mail:}
\verb"snasyrov@kpfu.ru"\vskip 0.5 cm

\end{document}